\documentclass[10pt]{amsart}
\usepackage{a4wide}
\usepackage{amssymb,amsmath,amsgen,amsxtra,dsfont,times} 
\begin{document}
%
%
%
\theoremstyle{definition}
\newtheorem{Definition}{Definition}[section]
\newtheorem{Example}[Definition]{Example}
\newtheorem{Examples}[Definition]{Examples}
\newtheorem{Remark}[Definition]{Remark}
\newtheorem{Remarks}[Definition]{Remarks}
\newtheorem{Caution}[Definition]{Caution}
\newtheorem{Conjecture}[Definition]{Conjecture}
\newtheorem{Question}[Definition]{Question}
\newtheorem{Questions}[Definition]{Questions}
\theoremstyle{plain}
\newtheorem{Theorem}[Definition]{Theorem}
\newtheorem{Proposition}[Definition]{Proposition}
\newtheorem{Lemma}[Definition]{Lemma}
\newtheorem{Corollary}[Definition]{Corollary}
\newtheorem{Fact}[Definition]{Fact}
\newtheorem{Facts}[Definition]{Facts}
\newtheoremstyle{voiditstyle}{3pt}{3pt}{\itshape}{\parindent}%
{\bfseries}{.}{ }{\thmnote{#3}}%
\theoremstyle{voiditstyle}
\newtheorem*{VoidItalic}{}
\newtheoremstyle{voidromstyle}{3pt}{3pt}{\rm}{\parindent}%
{\bfseries}{.}{ }{\thmnote{#3}}%
\theoremstyle{voidromstyle}
\newtheorem*{VoidRoman}{}
%
\newcommand{\prf}{\par\noindent{\sc Proof.}\quad}
\newcommand{\blowup}{\rule[-3mm]{0mm}{0mm}}
\newcommand{\cal}{\mathcal}
\newcommand{\Aff}{{\mathds{A}}}
\newcommand{\BB}{{\mathds{B}}}
\newcommand{\CC}{{\mathds{C}}}
\newcommand{\FF}{{\mathds{F}}}
\newcommand{\GG}{{\mathds{G}}}
\newcommand{\HH}{{\mathds{H}}}
\newcommand{\NN}{{\mathds{N}}}
\newcommand{\ZZ}{{\mathds{Z}}}
\newcommand{\PP}{{\mathds{P}}}
\newcommand{\QQ}{{\mathds{Q}}}
\newcommand{\RR}{{\mathds{R}}}
\newcommand{\Sphere}{{\mathds{S}}}
\newcommand{\lin}{\sim}
\newcommand{\num}{\equiv}
\newcommand{\dual}{\ast}
\newcommand{\iso}{\cong}
\newcommand{\homeo}{\approx}
\newcommand{\mm}{{\mathfrak m}}
\newcommand{\pp}{{\mathfrak p}}
\newcommand{\qq}{{\mathfrak q}}
\newcommand{\rr}{{\mathfrak r}}
\newcommand{\pP}{{\mathfrak P}}
\newcommand{\qQ}{{\mathfrak Q}}
\newcommand{\rR}{{\mathfrak R}}
%
%
\newcommand{\OO}{{\cal O}}
\newcommand{\into}{{\hookrightarrow}}
\newcommand{\onto}{{\twoheadrightarrow}}
\newcommand{\Spec}{{\rm Spec}\:}
\newcommand{\Pic}{{\rm Pic }}
\newcommand{\chit}{\chi_{\rm top}}
\newcommand{\KanDiv}{{\cal K}}
\newcommand{\perdef}{{\stackrel{{\rm def}}{=}}}
\newcommand{\Cycl}[1]{{\ZZ/{#1}\ZZ}}
\newcommand{\Sym}{{\mathfrak S}}
\newcommand{\EnSym}{{\cal S}}
\newcommand{\Alt}{{\mathfrak A}}
\newcommand{\Bra}{{\mathfrak B}}
\newcommand{\Free}{{\mathfrak F}}
\newcommand{\SymPro}{{\rm Sym}}
\newcommand{\ab}{{\rm ab}}
\newcommand{\Aut}{{\rm Aut}}
\newcommand{\Hom}{{\rm Hom}}
\newcommand{\Gal}{{\rm Gal}}
\newcommand{\Xgal}{{X_{{\rm gal}}}}
\newcommand{\Xaff}{{X^{{\rm aff}}}}
\newcommand{\Yaff}{{Y^{{\rm aff}}}}
\newcommand{\Xgalaff}{{X_{\rm gal}^{\rm aff}}}
\newcommand{\Xaffgal}{\Xgalaff}
\newcommand{\pitop}{{\pi_1}}
\newcommand{\Rgal}{{R_{\rm gal}}}
\newcommand{\fgal}{{f_{\rm gal}}}
\newcommand{\fgalast}{{f_{{\rm gal},\ast}}}
\newcommand{\Callg}{{C}}
\newcommand{\Caffine}{{C^{\rm aff}}}
\newcommand{\Cprojective}{{C^{\rm proj}}}
\newcommand{\barCaffine}{{{\bar{C}}^{\rm aff}}}
\newcommand{\barCprojective}{{{\bar{C}}^{\rm proj}}}
\newcommand{\myExt}{{\mathcal E}}
\newcommand{\myKernel}{{\mathcal K}}
\newcommand{\myBigExt}{\widetilde{{\mathcal E}}}
\newcommand{\myBigKernel}{\widetilde{{\mathcal K}}}
\newcommand{\catC}{{\mathcal C}}
\newcommand{\funF}{{\mathcal F}}
\newcommand{\Sets}{{\mathfrak{Sets}}}
\title[Fundamental Groups of Galois Closures]{Fundamental Groups of Galois Closures of\\Generic Projections}
\author{Christian Liedtke}
\address{Mathematisches Institut, Heinrich-Heine-Universit\"at, 40225
  D\"usseldorf}
\email{liedtke@math.uni-duesseldorf.de}
\date{November 2, 2005. {\it revised:} June 9, 2008}

\begin{abstract}
  For the Galois closure $\Xgal$ of a generic projection from
  a surface $X$, it is believed that $\pitop(\Xgal)$ 
  gives rise to new invariants of $X$.
  However, in all examples this group is surprisingly
  simple.
  In this article, we offer an explanation for this
  phenomenon:
  We compute a quotient of $\pitop(\Xgal)$ that depends on
  $\pitop(X)$ and data from the generic projection only.
  In all known examples, this quotient is in fact
  isomorphic to $\pitop(\Xgal)$.
  As a byproduct, we simplify the computations of
  Moishezon, Teicher and others.
\end{abstract}
\maketitle
%
\section*{Introduction}

One approach to a fine classification of algebraic surfaces
studies generic projections.
For that one embeds a given surface $X$ into some large projective 
space and projects it generically onto a plane.
If $D$ is the branch curve of such a projection, it
is possible to recover $X$ from the fundamental group 
$\pitop(\PP^2-D)$ and a little bit of extra data.
Even though these fundamental groups are extremely difficult to compute,
Moishezon \cite{mo} and later Teicher and others
developed tools (''braid monodromy factorisations``) to attack these 
computations systematically.
Eventually, one hopes to obtain new invariants of surfaces from these
groups.

One of these invariants of surfaces, 
which has itself a geometric interpretation, 
is the fundamental group of the corresponding Galois closure:
The Galois closure $\Xgal$ attached to a generic projection
from $X$ is again a smooth projective surface.
Its fundamental group $\pitop(\Xgal)$ is a subquotient of 
$\pitop(\PP^2-D)$ and therefore it is expected to contain
non-trivial information about the surface $X$.

Typically, $\Xgal$ is a surface of general type.
Miyaoka \cite{mi} proved that there are many surfaces with positive 
index among Galois closures of generic projections, i.e., their 
Chern numbers fulfil $c_1{}^2>2c_2$.
The  ''watershed-conjecture``, attributed
to Bogomolov, stated that surfaces of general type with 
positive index have infinite fundamental groups.
However, when Moishezon and Teicher \cite{mote} first computed 
$\pitop(\Xgal)$ for
generic projections from $\PP^1\times\PP^1$ they found 
simply connected surfaces of general type with positive index 
among them.
These were the first counter-examples to this conjecture.
Their computations are quite involved and use degeneration 
techniques, braid monodromy factorisations and lots of
combinatorial group theory.

Since then, the groups $\pitop(\Xgal)$ have been computed 
for many more examples.
For generic projections from minimal rational surfaces 
(except the Veronese surface $V_4$) the group $\pitop(\Xgal)$ 
is always finite Abelian.
Also, certain generic projections from $\PP^1\times E$, 
where $E$ is an elliptic curve, have been worked out.
Here, the groups $\pitop(\Xgal)$ are close to being free Abelian.
Why do these groups have such a simple structure?
\medskip

Here, we give an explaination:
Let $\Xgal$ be the Galois closure of a generic projection $f:X\to\PP^2$.
We let ${\cal L}:=f^\ast(\OO_{\PP^2}(1))$ be the very ample line bundle
that defines $f$.
We will assume that the degree $n:=\deg f$ is at least $5$.

In Theorem \ref{affinetheorem} and Theorem \ref{projectivemain} 
we determine the structure of 
$\pitop(\Xgal)/\Cprojective$, where $\Cprojective$ is a certain
subgroup generated by commutator and triple commutator relations.
Since our structure results need more notation and are a little
bit messy to state, we first give a couple of corollaries.

\begin{VoidItalic}[Corollary \ref{finitecor}]
  If $\pitop(X)$ is finite then so is $\pitop(\Xgal)/\Cprojective$.
\end{VoidItalic}

\begin{VoidItalic}[Corollary \ref{betti}]
  The rank of
  $H_1(\Xgal,\ZZ)/\overline{\Cprojective}$
  as an Abelian group is equal to $(n-1)$ times 
  the rank of $H_1(X,\ZZ)$.
\end{VoidItalic}

\begin{VoidItalic}[Proposition \ref{simplyconnected}]
  Let $X$ be a simply-connected surface.
  Let $d$ be the divisibility index of ${\cal L}$
  in $\Pic(X)$.
  Then there exists an isomorphism
  $$\begin{array}{llcl}
    &\pitop(\Xgal)/\Cprojective &\iso& (\Cycl{d})^{n-2} .
  \end{array}$$
\end{VoidItalic}

In all known examples, $\Cprojective$ is trivial.
These examples are generic projections from minimal
rational surface, which have been computed 
by Moishezon, Teicher and Robb
(\cite{mote}, \cite{mote2}, \cite{motero}), and
certain generic projections from $E\times\PP^1$, where $E$
is an elliptic curve, which have been computed by
Amram, Goldberg, Teicher and Vishne 
(\cite{ag}, \cite{agtv}).
We discuss these examples in Section \ref{examplesection}.
In particular, we compute $\pitop(\Xgal)$ modulo $\Cprojective$
without much effort: 
we obtain our quotients not in terms of generators and relations
but get a clear view on the structure of these groups.\medskip 

In view of Teicher's conjecture on virtual solvability of $\pitop(\Aff-D)$,
where $D$ is the branch curve of a generic projection from a surface
$X$ we contribute the following result which shows that there are
plenty of examples where $\pitop(\Aff^2-D)$ is not virtually solvable.

\begin{VoidItalic}[Proposition \ref{teicherconjecture}]
  If $\pitop(X)$ is not virtually solvable then neither is
  $\pitop(\Aff^2-D)$.
\end{VoidItalic}
\bigskip

To state the precise results, we have to introduce some notation.
We choose a generic line in $\PP^2$ and denote its complement by $\Aff^2$.
Then we consider the preimages of $\Aff^2$ on $X$ and $\Xgal$
and obtain the following diagram
$$\begin{array}{cccccc}
  \blowup
  \Xgal&\to&X&\stackrel{f}{\to}&\PP^2 \\
  \cup&&\cup&&\cup\\
   \Xaffgal&\to&\Xaff&\stackrel{f}{\to}&\Aff^2&.
\end{array}$$
For technical reasons (''monodromy at infinity``) it is easier to
compute the fundamental group $\pitop(\Xaffgal)$ first.

Modulo certain relations, the group
$\pitop(\Xaffgal)$ depends on $\pitop(\Xaff)$ and the degree 
$n$ of the generic projection only:
\begin{VoidItalic}[Theorem \ref{affinetheorem}]
  There exists a normal subgroup $\Caffine$ of $\pitop(\Xaffgal)$ and
  an isomorphism
  $$\begin{array}{ccc}
    \pitop(\Xaffgal)/\Caffine &\iso& \myBigKernel(\pitop(\Xaff),\,n) .
  \end{array}$$
  Here, $\myBigKernel(\pitop(\Xaff),n)$ 
  is a purely group theoretical construction, 
  which is an extension of $H_2(\pitop(\Xaff),\ZZ)$ 
  by a certain subgroup of $\pitop(\Xaff)^n$.
\end{VoidItalic}

The group $\pitop(\Xgal)$ is a quotient of $\pitop(\Xaffgal)$
by a cyclic subgroup.
The isomorphism of Theorem \ref{affinetheorem} is not well-behaved
with respect to forming this quotient.
In particular, if we denote by $\Cprojective$ the image of $\Caffine$
in $\pitop(\Xgal)$ then
$\pitop(\Xgal)/\Cprojective$ cannot be described
in terms of $n=\deg f$ and $\pitop(X)$ alone.

\begin{VoidItalic}[Theorem \ref{projectivemain}]
  There exists a cyclic group $Z$ and a central short exact
  sequence
  $$\begin{array}{ccccccccc}
      0 &\to &\displaystyle
      \frac{H_2(\pitop(\Xaff),\,\ZZ)}{Z}
      &\to& \pitop(\Xgal)/\Cprojective
      &\to& G &\to&1
    \end{array}$$
  where $G$ is a group that is itself a central extension
  $$\begin{array}{ccccccccc}
    0&\to&\ker\kappa_{n-1}
    &\to& G &\to&\myKernel(\pitop(X),\,n)
    &\to&1.
  \end{array}$$
  Here, $\myKernel(\pitop(X),n)$ is a certain
  subgroup of $\pitop(X)^n$.
  The group $\ker\kappa_{n-1}$ is a subgroup of $K^{n-1}$, where
  $K$ is a cyclic group.
\end{VoidItalic}

We show how to compute $\ker\kappa_{n-1}$ and 
$H_2(\pitop(\Xaff),\ZZ)$ in practice in 
Section \ref{examplesection}.
Proposition \ref{simplyconnected}, Corollary \ref{betti}
and Corollary \ref{finitecor} quoted above are applications 
of this theorem.

However, the main question that remains is the following

\begin{VoidRoman}[Question \ref{mainquestion}]
  Is $\Caffine$ trivial for every generic projection
  of degree $n\geq5$?
\end{VoidRoman}

This is true in all known examples.
However, since not many examples with nontrivial 
$\pitop(X)$ have been computed, we do not dare to call it 
a conjecture.
The assumption on the degree is necessary:
a generic projection from the Veronese surface $V_4$ in 
$\PP^5$, which has only degree $4$, leads to a nontrivial
$\Caffine$, cf. Remark \ref{veronese}.
\medskip

We now relate Question \ref{mainquestion} to the computations of
Moishezon, Teicher and others.
To calculate $\pitop(\Xaffgal)$,
one first needs a presentation of the group $\pitop(\Aff^2-D)$, where
$D$ is the branch curve of the generic projection.
Using degeneration techniques,
one then shows that certain commutator and triple commutator relations
hold between the standard generators of $\pitop(\Aff^2-D)$, 
cf. \cite[Proposition 2]{mote}.
In our setup this amounts to proving that the group
$\Caffine$ is trivial.
Hence answering Question \ref{mainquestion} is a crucial step
in the calculation of $\pitop(\Xaffgal)$.
Once we know that $\Caffine$ is trivial, we can determine 
$\pitop(\Xgal)$, which already simplifies the existing computations
a lot.
\medskip

In any case, our results shed new light on $\pitop(\Xgal)$ 
and open new perspectives.
\begin{itemize}
\item 
  Theorem \ref{projectivemain} and its corollaries 
  tell us that 
  $\pitop(\Xgal)/\Cprojective$ is not too far away 
  from $\pitop(X)$.
  For example, if $X$ is simply connected, 
  $\pitop(\Xgal)/\Cprojective$ 
  does not contain interesting information on $X$.
  Even if this group is hard to determine,
  our structure results and Question \ref{mainquestion}
  suggest that $\pitop(\Xgal)$ does not give rise to
  new invariants of $X$.
\item
  The group $H_2(\pitop(\Xaff),\ZZ)$ maps to 
  the commutator subgroup of $\pitop(\Xgal)/\Cprojective$ 
  (Theorem \ref{projectivemain}).
  Unfortunately, the known examples are too special so 
  that this homological contribution has not yet appeared.
  Since it is still not known, what kind of infinite groups 
  can occur as fundamental groups of surfaces, one could use
  iterated Galois closures of generic projections and 
  Corollary \ref{nilpotentcor} to construct surfaces 
  with nilpotent fundamental groups of 
  arbitrary large nilpotency class.
\end{itemize}
\medskip

The article is organised as follows.

In Section \ref{GenProj} we review some facts about 
generic projections and their Galois closures.

In Section \ref{TwoGroups} we introduce two purely group theoretical
constructions.
Both depend on a group $G$ and a natural number $n\geq3$.
The first construction, called $\myKernel(G,n)$ is a certain 
subgroup of $G^n$ and quite easy to calculate.
The second one, called $\myBigKernel(G,n)$ is a central extension of
$H_2(G,\ZZ)$ by $\myKernel(G,n)$ and quite difficult to determine in
general.
There is a surprising connection between $\myBigKernel(G,n)$ and the
theory of central extensions whose meaning is unclear at the moment.
We refer to \cite{semiversal} for purely group theoretical properties
of $\myBigKernel(G,n)$.

In Section \ref{Mainaffine} we consider the short exact sequence 
\begin{equation}{\tag{\ref{moteeq}}}
\begin{array}{ccccccccc}
  1&\to&\pitop(\Xaffgal)&\to&\pitop(\Aff^2-D)/\ll\Gamma_i{}^2\gg
  &\to&\Sym_n&\to&1,
\end{array}
\end{equation}
which is also the starting point of the calculations 
of Moishezon and Teicher in \cite{mote}.
Here, $D$ is the branch curve of the generic projection
$f:X\to\PP^2$, the $\Gamma_i$'s are standard generators for 
$\pitop(\Aff^2-D)$ and $\Sym_n$ is the symmetric group on
$n:=\deg f$ letters.
We introduce a normal subgroup $\Caffine$ of
$\pitop(\Aff^2-D)$.
This group is trivial if certain commutator and triple
commutator relations hold between the $\Gamma_i$'s, cf.
Definition \ref{caffinedef}.
After taking the quotient of
the short exact sequence (\ref{moteeq}) by $\Caffine$, the
resulting short exact sequence splits in a very nice way.
Using the algorithm of Zariski and van~Kampen and the results
of Section \ref{TwoGroups}, we obtain Theorem \ref{affinetheorem}.

In Section \ref{Mainprojective} we determine the structure of
$\pitop(\Xgal)/\Cprojective$ and give some applications.

In Section \ref{examplesection} we compute these quotients 
of $\pitop(\Xgal)$ and $\pitop(\Xaffgal)$ for all examples
that are known to the author.
In particular, we demonstrate how to avoid combinatorial
group theory once our theorems are at disposal.

\begin{VoidRoman}[Acknowledgements]
This article extends the main results of my
Ph.D. thesis \cite{phd}.
I thank my supervisor Gerd~Faltings for 
suggesting this interesting topic and the discussions on it.
Also, I thank the  
Max-Planck-Institut in Bonn for hospitality and financial 
support.
Finally, I thank the referee for pointing out a couple
of inaccuracies.
\end{VoidRoman}

\section{Generic projections and their Galois closures}
\label{GenProj}

In this section we recall some general facts about generic 
projections and their Galois closures.
Most of these results are well-known.
We follow \cite[Section 2]{fa}.
\begin{Definition}
  Let $X$ be a smooth complex projective surface.
  A line bundle ${\cal L}$ on $X$ is called
  {\em sufficiently ample} if
  \begin{enumerate}
    \item ${\cal L}$ is an ample line bundle with self-intersection
      number at least $5$.
    \item For every closed point $x\in X$ the global sections of
      $\cal L$ generate the fibre 
      $${\cal L}_x/\mm_x{}^4{\cal L}.$$
    \item For any pair $\{x,y\}$ of distinct closed points of $X$
      the global sections of $\cal L$ generate the direct sum
      $${\cal L}_x/\mm_x{}^3{\cal L}\oplus
      {\cal L}_y/\mm_y{}^3{\cal L}.$$
    \item For any triple $\{x,y,z\}$ of distinct closed points
      the global sections of $\cal L$ generate the direct sum
      $${\cal L}_x/\mm_x{}^2{\cal L}\oplus
      {\cal L}_y/\mm_y{}^2{\cal L}\oplus
      {\cal L}_z/\mm_z{}^2{\cal L}.$$
  \end{enumerate}
\end{Definition}
The following remark was already made in \cite{fa}.
\begin{Remark}
\label{suffample}
 The tensor product of five very ample line bundles is
 sufficiently ample.
 \end{Remark}
\begin{Definition}
  \label{genericdef}
  Let $\cal L$ be a sufficiently ample line bundle on the 
  surface $X$.
  We call a three-dimensional linear subspace 
  $E\subseteq H^0(X,{\cal L})$ {\em generic} if
  \begin{enumerate}
  \item $E$ generates $\cal L$, i.e., the associated rational
    map is in fact a finite morphism
    $$\begin{array}{ccccc}
      f=f_E&:&X&\to&\PP(E)\iso\PP^2
    \end{array}$$
    of degree $n$ equal to the self-intersection number of $\cal L$.
  \item The ramification locus $R$ of $f_E$ is a smooth and
    ample curve on $X$.
    The ramification index at a generic point is $2$.
  \item The branch locus  $D$ of $f_E$ is a curve on 
    $\PP^2$ with at worst cusps and nodes as singularities.
  \item The restriction $f_E|_R:R\to D$ is birational.
  \end{enumerate}
  We call the finite morphism $f_E$ associated to a generic subspace
  of $H^0(X,{\cal L})$ a {\em generic projection}.
\end{Definition}

Let $\GG(k,V)$ be the Grassmannian parametrising $k$-dimensional
linear subspaces of the vector space $V$.
To justify the name in the previous definition we have

\begin{Proposition}
  Let $\cal L$ be a sufficiently ample line bundle on 
  the surface $X$.
  Then there exists a dense and open subset 
  $G'$ of $\GG(3,H^0(X,{\cal L}))$
  such that all $E\in G'$ are generic in the sense of
  Definition \ref{genericdef}.
\end{Proposition}

\prf
A proof is given in \cite[Proposition 1]{fa}.
For the conclusions there it is necessary that the line
bundle ${\cal K}_X\otimes{\cal L}^{\otimes 3}$ is ample.
In our setup, this follows from the assumptions on the 
self-intersection number of $\cal L$ and Reider's 
theorem \cite{re}.
\qed\medskip

Let $f:X\to P$ be a finite morphism between normal complex 
surfaces.
Let $L$ be the Galois closure of the field extension of
function fields $K(P)\,\into\, K(X)$ induced by $f$. 
Let $\Xgal$ be the normalisation of $X$ or, equivalently, of
$P$ in $L$.
\begin{Definition}
  Given a finite morphism $f:X\to P$ of normal complex surfaces
  we call the normal surface $\Xgal$ together with the 
  induced morphism $\fgal:\Xgal\to P$ the {\em Galois closure}
  of $f$.
\end{Definition}
The Galois closure of a field extension $K/k$ of degree $n$ 
is isomorphic to a quotient of the $n$-fold tensor product of
$K$ with itself over $k$.
Similarly (cf. \cite[Expos\'e V.4g]{sga1}), the Galois closure of a
generic projection is isomorphic to the closure of a certain subset
of $X^n$:
$$\begin{array}{ccccc}
  \Xgal &\iso& \overline{\{ (x_1,...,x_n)\,|\, x_i\neq x_j, f(x_i)=f(x_j), 
  \forall i\neq j\}}
  &\subseteq& X^n\,.
\end{array}$$

For Galois closures of generic projections we have the
following proposition, which is proved for our setup
in \cite[Proposition 1]{fa}.
\begin{Proposition}
  \label{goodsubspace}
  Let $\cal L$ be a sufficiently ample line bundle on the
  surface $X$ with self-intersection $n$.
  Then there exists an open and dense subset $G$ of
  $\GG(3,H^0(X,{\cal L}))$ such that for all $E\in G$
  \begin{enumerate}
  \item The map $f_E$ associated to $E$ is a generic projection
    of degree $n$.
  \item The Galois closure $\Xgal$ of $f_E$ is a smooth projective
    surface.
  \item The Galois group $\Gal(K(\Xgal)/K(\PP^2))$ is the symmetric
    group $\Sym_n$. 
  \item The branch curve $D$ of $f_E$ is an irreducible curve
    in $\PP^2$.
  \end{enumerate}
\end{Proposition}

\begin{Definition}
  We call a map $f_E$ as in Proposition \ref{goodsubspace}
  a {\em good generic projection} from $X$.
\end{Definition}

Explicit formulae for the degree of the branch curve $D$ and the number of
cusps and nodes on $D$ can be found in \cite[Section 4]{fa} or
\cite[Chapter 0]{mote}.

\section{Two group theoretical constructions}
\label{TwoGroups}

This section contains the group theoretical heart of this article.
Given a group $G$ and a natural number $n\geq3$ we introduce the two
associated groups $\myKernel(G,n)$ and $\myBigKernel(G,n)$.
The second construction is related to the theory of central extensions.
We refer to \cite{semiversal} for details, further properties and examples.

\subsection*{The first construction}

We let $G$ be a group and $n\geq3$ be a natural number.
We denote by $G^\ab:=G/ [G,G]$ the Abelianisation of $G$.
Then we define a map
$$\begin{array}{ccccc}
\psi&:&G^n&\to&G^\ab\\
&&(g_1,...,g_n)&\mapsto&\overline{g_1\cdot...\cdot g_n}
\end{array}$$
which is a homomorphism since $G^\ab$ is an Abelian group.

\begin{Definition}
  For a group $G$ and a natural number $n\geq3$ we define
  $\myKernel(G,\,n)$ to be the kernel of the homomorphism
  $\psi:G^n\to G^\ab$.
\end{Definition}

For Abelian groups it is particularly easy to determine $\myKernel(G,n)$.
\begin{Proposition}
  \label{mykernelAbelian}
  There exists a natural isomorphism
  $\myKernel(G,n)^\ab\iso\myKernel(G^\ab,n)$.

  If $G$ is Abelian then
  $\myKernel(G,n)\iso G^{n-1}$.
  This isomorphism is not canonical.
\end{Proposition}
\prf
The surjection $G\to G^\ab$ and the universal property of Abelianisation
imply that there is a natural surjective homomorphism
$\myKernel(G,n)^\ab\to\myKernel(G^\ab,n)$.
An element of the kernel $\myKernel(G,n)\to\myKernel(G^\ab,n)$ is also
an element of the kernel $G^n\to(G^\ab)^n$, which is $[G,G]^n$.
Using $n\geq3$ we may write
\begin{equation}
\label{commutatoreq}
\begin{array}{ccc}
  ([h_1,h_2],1,...,1) &=& [ (h_1,h_1{}^{-1},1,...,1),\,(h_2,1,h_2{}^{-1},...,1) ] .
\end{array}
\end{equation}
Thus $[G,G]^n$ is not only a subgroup of $\myKernel(G,n)$ but also lies
inside the commutator subgroup of $\myKernel(G,n)$.
Hence of the kernel $\myKernel(G,n)\to\myKernel(G^\ab,n)$ is the commutator
subgroup of $\myKernel(G,n)$ and we are done.

Let $G$ be Abelian.
Then the map
$$\begin{array}{ccc}
  G^{n-1}&\to& G^n\\
  (g_1,...,g_{n-1})&\mapsto&(g_1,...,g_{n-1},(g_1\cdot...\cdot g_{n-1})^{-1})
\end{array}$$
defines a homomorphism.
It is injective with image $\myKernel(G,n)$.
\qed\medskip

The symmetric group $\Sym_n$ acts on $G^n$ by permutation
of the $n$ factors.
This action induces an action $\varphi$ 
of $\Sym_n$ on $\myKernel(G,n)$.
With respect to this action we form the semidirect product
$$\begin{array}{ccc}
  \myExt(G,n) &:=& \myKernel(G,n)\rtimes_\varphi\Sym_n\,.
\end{array}$$
For $\vec{g}\in\myKernel(G,n)$ and
$\sigma\in\Sym_n$ we will write $[\vec{g},\sigma]$ for
$\vec{g}\cdot \varphi(\sigma)(\vec{g}^{-1})$.

\subsection*{The main construction}

Again, we let $G$ be a group and $n\geq3$ be a natural number.
We choose a presentation $G\iso F/N$, where $F$ is a free group.
Then $\myKernel(N,n)$ is a subgroup of $\myKernel(F,n)$, which is
a subgroup of $F^n$.

We denote by $\ll\myKernel(N,n)\gg$ the subgroup normally generated
by $\myKernel(N,n)$ inside $F^n$.
Using formula (\ref{commutatoreq}), it is not hard to verify that the
group $\ll\myKernel(N,n)\gg$ is equal to the subgroup normally
generated by $\myKernel(N,n)$ inside $\myKernel(F,n)$.

\begin{Definition}
  Let $G$ be a group and $n\geq3$ be a natural number.
  We define
  $$\begin{array}{ccc}
    \myBigKernel(G,n) &:=& \myKernel(F,n)/\ll\myKernel(N,n)\gg .
  \end{array}$$
\end{Definition}

\begin{Theorem}
  \label{mybigkernelthm}
  The group $\myBigKernel(G,n)$ does not depend upon the choice
  of a presentation.
  There exists a central short exact sequence
  \begin{equation}
  \label{mybigkerneleq}
  \begin{array}{ccccccccc}
    0&\to&H_2(G,\,\ZZ)&\to&\myBigKernel(G,\,n)&\to&\myKernel(G,n)&\to&1,
  \end{array}
  \end{equation}
  where $H_2(G,\ZZ)$ denotes the second group homology with
  integer coefficients.
  The group $H_2(G,\ZZ)$ lies inside the commutator
  subgroup of $\myBigKernel(G,n)$.
\end{Theorem}
\prf
We choose a presentation $G\iso F/N$ and abbreviate the normal
closure of $\myKernel(N,n)$ in $F^n$ by $R$.
Also, we will simply write $H_2(G)$ for $H_2(G,\ZZ)$.

First, we prove the short exact sequence (\ref{mybigkerneleq}):
We denote by $\pi$ the projection of $F^n$ onto its last
$n-1$ factors. 
By abuse of notation we will denote its restriction to 
$\myKernel(F,n)$ also by $\pi$.
Then we have a short exact sequence
$$\begin{array}{ccccccccc}
  1&\to&[F,F]&\to&\myKernel(F,n)&\stackrel{\pi}{\to}&F^{n-1}&\to&1.
\end{array}$$
An easy computation with commutators shows that $R\cap\ker\pi=[F,N]$.
Via $\pi$ we obtain the following diagram with exact rows and columns:
\begin{equation}
\label{mainprfeq}
  \begin{array}{ccccccccc}
  1&\to&[F,N]&\to&R&\to&N^{n-1}&\to&1 \\
  &&\downarrow&&\downarrow&&\downarrow\\
  1&\to&N\cap[F,F]&\to&N^n\cap\myKernel(F,n)&\to&N^{n-1}&\to&1 \\
  &&\downarrow&&\downarrow&&\downarrow\\
  1&\to&[F,F]&\to&\myKernel(F,n)&\stackrel{\pi}{\to}&F^{n-1}&\to&1 
\end{array}
\end{equation}
Taking quotients of successive rows we exhibit the group $\myKernel(F,n)/R$ as an
extension of $(N\cap[F,F])/[F,N]$ by $\myKernel(F,n)/(N^n\cap\myKernel(F,n))$.
The latter group is isomorphic to $\myKernel(G,n)$.
By Hopf's theorem \cite[Theorem II.5.3]{br} the group
$(N\cap [F,F])/ [F,N]$ is isomorphic to $H_2(G)$.
Hence we obtain an extension
$$\begin{array}{ccccccccc}
  1&\to&H_2(G)&\to&\myKernel(F,n)/R&\to&\myKernel(G,n)&\to&1.
\end{array}$$

Now we will show that this extension is central:
Every element of $H_2(G)$ can be lifted to an element $x\in\myKernel(F,n)$
of the form $\vec{x}:=(x,1,...,1)$ with $x\in N\cap[F,F]$.
For $\vec{y}:=(y_1,...,y_n)\in\myKernel(F,n)$ we compute
$$
\vec{y}\,\vec{x}\,\vec{y}^{-1}\,=\,
(\underbrace{[y_1,x]}_{\in[F,N]},1,...,1)\cdot(x,1,...,1) \,\equiv\,
\vec{x}\,\mod[F,N].
$$
Hence $H_2(G)$ lies inside the centre of $\myKernel(F,n)/R$.

We now prove that $\myBigKernel(G,n)$ is well-defined:
Let $\alpha:F/N\iso F'/N'$ be another presentation for $G$.
We may lift this isomorphism to a map $\varphi:F\to F'$.
Then $\varphi$ maps $N$ to $N'$ and hence $\myKernel(N,n)$ to
$\myKernel(N',n)$.
Let $R'$ be the normal closure of $\myKernel(N',n)$ inside $F'^n$.
Then $\varphi$ induces a homomorphism
$$\begin{array}{ccccc}
\overline{\varphi} &:&\myKernel(F,n)/R &\to&\myKernel(F',n)/R'
\end{array}$$
Let $\varphi'$ be another map lifting $\alpha$ to a homomorphism
from $F$ to $F'$.

Elements of the form $(f,f^{-1},1,....,1)$ generate $\myKernel(F,n)$.
Hence it is enough to compare the induced maps on these elements.
For $f\in F$ there exists $n_f'\in N'$ such that
$\varphi(f)=\varphi'(f) n_f'$.
Hence
$$\begin{array}{lcl}
  \varphi((f,f^{-1},...)) &=& 
  (\varphi'(f)n_f', n_f'^{-1}\varphi'(f)^{-1},...) \\
  &=&\varphi'((f,f^{-1},...)) 
  \underbrace{(n_f',\varphi'(f)n_f'^{-1}\varphi'(f)^{-1},...)}_{\in R'}.
\end{array}$$
Hence the induced maps coincide.

In particular, if we choose $F=F'$ and $N=N'$ with $\alpha$ and
$\varphi$ the identity then every other lift $\varphi'$ of the
identity induces the identity on $\myKernel(F,n)/R$.

Coming back to the general case, we let $F/N$ and $F'/N'$ again be
two presentations of $G$ and let $\alpha$ be an isomorphism between
them.
Then $\alpha$ and $\alpha^{-1}$ induce maps between $\myKernel(F,n)/R$
and $\myKernel(F',n)/R'$ such that the composites of these induced maps
are the identity by the previous paragraph.
Hence $\alpha$ induces an isomorphism from $\myKernel(F,n)/R$ to
$\myKernel(F',n)/R'$.
Thus, $\myBigKernel(G,n)$ is well-defined.

Taking the quotient of the top row by the bottom row of (\ref{mainprfeq})
we obtain an exact sequence
$$\begin{array}{ccccccccc}
  1&\to&[F,F]/[F,N]&\to&\myKernel(F,n)/R&\to&G^{n-1}&\to&1.
\end{array}$$
The inclusion of $H_2(G)$ into $\myKernel(F,n)/R$ factors over 
$[F,F]/[F,N]$.
The group $[F,F]$ lies inside the commutator
subgroup of $\myKernel(F,n)$, cf. formula (\ref{commutatoreq}).
Hence $H_2(G)$ lies inside the commutator subgroup of $\myKernel(F,n)/R$.
\qed\medskip

Whenever the symmetric group $\Sym_n$ acts on a group $X$ via some
homomorphism $\varphi:\Sym_n\to\Aut(X)$ we may form the semidirect
product $X\rtimes\Sym_n$ with respect to this action.
For $x\in X$ and $\sigma\in\Sym_n$ the commutator $[x,\sigma]$ in
$X\rtimes\Sym_n$ is equal to 
$x\sigma x^{-1}\sigma^{-1}=x\cdot\varphi(\sigma)(x^{-1})$ by 
definition of the semidirect product.
For a subgroup $S$ of $\Sym_n$ we define $[X,S]$ to be the
subgroup of $X\rtimes\Sym_n$ generated by all elements $[x,s]$
with $x\in X$ and $s\in S$.
It is easy to see that $[X,S]$ is a normal subgroup of $X$.

We already mentioned that there is an action of $\Sym_n$ on
$\myKernel(G,n)$ for every group $G$.
In the following we embed $\Sym_{n-1}$ into $\Sym_n$ as the
subgroup of those permutations that fix the first letter.

\begin{Corollary}
  \label{recover}
  There exists a $\Sym_n$-action on $\myBigKernel(G,n)$ such 
  that the homomorphism from $\myBigKernel(G,n)$ onto 
  $\myKernel(G,n)$ becomes equivariant.
  There exists an isomorphism
  $$\begin{array}{ccc}
    \myBigKernel(G,\,n)\,/\,%
    [\myBigKernel(G,\,n),\,\Sym_{n-1}] &\iso& G .
  \end{array}$$ 
\end{Corollary}
\prf
Given a presentation $G\iso F/N$ it is easy to 
check that the $\Sym_n$-action on $\myKernel(F,n)$ induces
an action on $\myBigKernel(G,n)$ that does not depend on the
choice of the presentation.

We view $\myKernel(F,n)$ as a subgroup of $F^n$.
We let $\Sym_{n-1}$ be the subgroup of
$\Sym_n$ of those permutations that fix the first letter. 
The group $[\myKernel(F,n),\Sym_{n-1}]$
is a normal subgroup of $\myKernel(F,n)$ consisting of
those elements with a trivial entry in the first coordinate.
Hence
$$\begin{array}{ccc}
  \myKernel(F,n)\,/\,\left[\myKernel(F,\,n),\,\Sym_{n-1}\right]
  &\iso& F .
\end{array}$$
Of course, the same holds true if we replace $F$ by $N$.
We have isomorphisms
$$\frac{\myBigKernel(G,\,n)}{\left[\myBigKernel(G,\,n),\,\Sym_{n-1}\right]}
  \,\iso\,
  \frac{\myKernel(F,\,n)}{\left[\myKernel(F,\,n),\,\Sym_{n-1}\right]}\,/\,
  \frac{\myKernel(N,\,n)}{\left[\myKernel(N,\,n),\,\Sym_{n-1}\right]}
  \,\iso\,
  F/N\,\iso\,G .$$
This proves our assertion.
\qed

\begin{Corollary}
  \label{bigcyclic}
  If $G$ is cyclic then
  $\myBigKernel(G,n)\iso G^{n-1}$.
  This isomorphism is not canonical.
\end{Corollary}
\prf
If $G$ is cyclic then $H_2(G,\ZZ)$ vanishes.
Hence $\myBigKernel(G,n)$ is isomorphic to $\myKernel(G,n)$,
which is isomorphic to $G^{n-1}$ by Proposition \ref{mykernelAbelian}.
\qed

\begin{Corollary}
  \label{bigabelianise}
  There exist isomorphisms
  $\myBigKernel(G,n)^\ab\iso\myKernel(G,n)^\ab\iso (G^\ab)^{n-1}$.
\end{Corollary}
\prf
The first isomorphism follows from the fact that
$H_2(G,\ZZ)$ lies inside the commutator subgroup of
$\myBigKernel(G,n)$.
The second isomorphism follows from Proposition \ref{mykernelAbelian}.
\qed

\subsection*{A presentation of $\myExt(F,n)$ and affine subgroups}

Let $F$ be a free group.
Then $\myKernel(F,n)$ is a subgroup of $F^n$.
We recall that there exists an action $\varphi$ of the
symmetric group $\Sym_n$ on $\myKernel(F,n)$.
For $\vec{f}\in\myKernel(F,n)$ and $\sigma\in\Sym_n$, 
we refer to $\varphi(\sigma)(\vec{f})$ as the element
that is $\sigma$-conjugate to $\vec{f}$.

\begin{Definition}
  \label{affinesubdef}
  Let $F$ be a free group.
  A normal subgroup of $\myExt(F,n)$ group is called 
  {\em affine subgroup} if it is normally generated
  by elements of the form
  $(r,r^{-1},1,...,1)$ and their $\Sym_n$-conjugates.
\end{Definition}

Clearly, affine subgroups of $\myKernel(F,n)$ are $\Sym_n$-invariant.
The quotient of $\myKernel(F,n)$ by an affine subgroup has a
particularly nice structure:

\begin{Theorem}
  \label{affinesubthm}
  Let $F$ be a free group and let $R$ be an affine subgroup of 
  $\myKernel(F,n)$.
  Let $p:\myKernel(F,n)\to F$ be the projection onto the first
  factor.
  We define $N:=p(R)$.
  Then there exists an isomorphism
  $$\begin{array}{ccc}
    \myKernel(F,\,n)/R &\iso& \myBigKernel(F/N,\,n).
  \end{array}$$
  In particular, this quotient is completely determined by
  $F/N$ and $n$.
\end{Theorem}
\prf
As $p$ is surjective, $N=p(R)$ is a normal subgroup 
of $F$.
Using $n\geq3$ and formula (\ref{commutatoreq}) it is easy
to see that $\myKernel(N,n)$ is contained in $R$.
As $R$ is a normal subgroup of $\myKernel(F,n)$,
also the normal closure $\ll\myKernel(N,n)\gg$
is contained in $R$.

Conversely, $R$ is normally generated by elements of the form
$(r,r^{-1},...)$ and their $\Sym_n$-conjugates.
But if $(r,r^{-1},1...)$ lies in $R$ then $r$ lies in $N$.
We conclude that $R$ is contained in $\ll\myKernel(N,n)\gg$.
Thus $R$ is equal to $\ll\myKernel(N,n)\gg$ and so
$\myKernel(F,n)/R$ is isomorphic to $\myBigKernel(F/N,n)$.
\qed\medskip

Let $F_{d-1}$ be the free group of rank $d-1$.
We finish this section by giving a presentation of $\myExt(F_{d-1},n)$:
Let $S_d$ be the free group of rank $d$ with basis $s_i$, $i=1,...,d$.
We form the free product of the symmetric group $S_d$ with $\Sym_n$.
We define
$$\begin{array}{cccccl}
  \psi &:& S_d\ast\Sym_n &\to& \Sym_n \\
  &&s_i&\mapsto&(1\,2)&i=1,...,d\\
  &&\sigma&\mapsto&\sigma&\forall\sigma\in\Sym_n .
\end{array}$$
We denote by
$[a,b]:=aba^{-1}b^{-1}$ and
$\langle a,b\rangle:=abab^{-1}a^{-1}b^{-1}$ 
commutators and triple commutators, respectively.
Then we let $R_{d,n}$ be the subgroup of $S_d\ast\Sym_n$
normally generated by the elements
$$\begin{array}{ll}
  s_i{}^2 &\mbox{ for $i=1,...,d$ }\\
  s_1\cdot (1\,2)^{-1} & 
  \mbox{i.e., we identify $s_1$ with $(1\,2)$}\\
  
  \left[s_i,\,\tau\right] & \mbox {for transpositions $\tau\in\Sym_n$ 
    disjoint from $(1\,2)$}\\
  
  \langle s_i,\,\tau\rangle & \mbox {for transpositions $\tau\in\Sym_n$
    having one index}\\
  &\mbox{in common with $(1\,2)$}\\
  
  \left[\sigma s_i\sigma^{-1},\,s_j\right] & 
  \mbox {if $\psi(\sigma s_i\sigma^{-1})$ and $\psi(s_j)$ are
    disjoint}\\
  \langle\sigma s_i\sigma^{-1},\,s_j\rangle & 
  \mbox {if $\psi(\sigma s_i\sigma^{-1})$ and $\psi(s_j)$ have one
    index in common}
\end{array}$$

\begin{Definition}
  \label{ensymdef}
  We define $\EnSym_n(d)\,:=\,(S_d\ast\Sym_n)/R_{n,d}$.
\end{Definition}
It is easy to see that the homomorphism $\psi$ factors over
the quotient by $R_{n,d}$.
By abuse of notation we will call the induced map on $\EnSym_n(d)$ 
again $\psi$.
We may think of $\EnSym_n(d)$ as a symmetric group with
$d$ ''layers``. \medskip

Let $F_{d-1}$ be the free group of rank $d-1$ with basis
$f_2,...,f_d$ (mind the numbering!).
The group $\Sym_n$ acts on $(F_{d-1})^n$ by permuting the factors.
We denote the corresponding semidirect product by
$(F_{d-1})^n\rtimes\Sym_n$.
It is easy to check that the map
$$\begin{array}{cccccl}
  \varphi&:&\EnSym_n(d)&\to&(F_{d-1})^n\rtimes\Sym_n\\
  &&s_1&\mapsto&(1\,2)\\
  &&s_i&\mapsto&(f_i,f_i{}^{-1},1,...,1)\cdot(1\,2)&i=2,...,d\\
  &&\sigma&\mapsto&\sigma&\forall \sigma\in\Sym_n
\end{array}$$
defines a homomorphism of groups.
Its image lies inside $\myExt(F_{d-1},n)$.
Moreover, this latter group is generated
by the $\Sym_n$-conjugates of the elements $(f_i,f_i{}^{-1},1,...,1)$.
Hence $\varphi$ maps surjectively onto
$\myExt(F_{d-1},n)$.
But even more is true

\begin{Theorem}
  \label{rtv}
  For $n\geq5$ the map $\varphi$ defines an isomorphism
  $$\begin{array}{ccc}
    \EnSym_n(d) &\iso& \myExt(F_{d-1},\,n).
  \end{array}$$
  In particular, $\EnSym_n(d)$ defines a presentation
  of $\myExt(F_{d-1},n)$.
\end{Theorem}
\prf
Our original proof applied a Reidemeister-Schreier 
rewriting process to the kernel $\ker(\psi)$ of $\EnSym_n(d)$. 
We obtained an infinite set of relations and 
checked that all of them follow from relations 
in $(F_{d-1})^n$.
This then implies that $\varphi$ is injective.
The computations are straight forward, but quite lengthy.

In the meantime, the article \cite{rtv} appeared independently
from our work.
So instead of giving our proof, we refer to
\cite[Theorem  5.7]{rtv} instead.
\qed

\section{The fundamental group of the affine piece}
\label{Mainaffine}

Let $f:X\to\PP^2$ be a good generic projection of degree $n$
given by a sufficiently ample line bundle as defined
in Section \ref{GenProj}.
We denote by $\fgal:\Xgal\to\PP^2$ its Galois closure.
We choose a generic line $\tilde{\ell}$ in $\PP^2$, i.e., 
a line that intersects the branch curve $D$ of $f$ in 
$\PP^2$ transversely.

We define $\Xaff:=X-f^{-1}(\tilde{\ell})$ and
$\Xaffgal:=\Xgal-\fgal^{-1}(\tilde{\ell})$.
Then we obtain the following morphisms
$$\begin{array}{lcccccc}
  \blowup&\Xgal &\to &X&\stackrel{f}{\to}&\PP^2\\
  &\cup&&\cup&&\cup\\
  &\Xaffgal &\to &\Xaff&\stackrel{f}{\to}&\Aff^2 &.
\end{array}$$

We choose a second line $\ell$ in $\Aff^2:=\PP^2-\tilde{\ell}$
that intersects $D$ and $\tilde{\ell}$ transversely.
Then we choose a system of simple loops $\Gamma_i$ around the
intersection points $\ell\cap D$ inside $\ell-D$.
There are $d:=\deg D$ such simple loops.
The inclusion map of topological spaces induces a surjective homomorphism
$$\begin{array}{ccc}
\pitop(\ell-(D\cap\ell))&\onto&\pitop(\Aff^2-D).
\end{array}$$
The fundamental group on the left is a free group of rank $d$ with
basis $\left\{\Gamma_i\right\}$, $i=1,...,d$.

By definition of a generic projection (Definition \ref{genericdef}),
the branch curve $D$ has at worst cusps and nodes
as singularities. 

\begin{Lemma}
  \label{vankampenlemma}
  The simple loops $\Gamma_i$, $i=1,...,\deg D$ generate
  $\pitop(\Aff^2-D)$.
  All relations in this group follow from relations
  of the following type
  $$
  \gamma\Gamma_i\gamma^{-1}\cdot\Gamma_j{}^{-1},\,\,
  \left[\gamma\Gamma_i\gamma^{-1},\,\Gamma_j\right],\,\,
  \langle\gamma\Gamma_i\gamma^{-1},\,\Gamma_j\rangle
  $$
  for elements $\gamma\in\pitop(\Aff^2-D)$.
\end{Lemma}
\prf
This follows from the theorem of Zariski and van~Kampen
on fundamental groups of complements. 
A modern proof is given in \cite[Th\'eor\`eme 4.1]{ch}.
The relations occurring if the curve has at worst cusps and 
nodes are determined in \cite[Partie 6.2]{ch}.
\qed\medskip

We lift loops around $x_0$ to points of $f^{-1}(x_0)$
and obtain paths in $\Xaff$.
Clearly, the endpoints of these paths are again points of
that fibre.
Thus every loop around $x_0$ defines a
permutation of $f^{-1}(x_0)$, which is a set with $n$ elements.
So there exists a homomorphism $\psi$ from
$\pitop(\Aff^2-D)$ to the symmetric group $\Sym_n$.

The Galois group of the extension of the function
fields of $\Xgal$ over $\PP^2$ is isomorphic to
$\Sym_n$ by Proposition \ref{goodsubspace}. 
It follows that 
$$\begin{array}{ccccc}
  \psi &:&\pitop(\Aff^2-D,\,x_0)&\onto&\Sym_n
\end{array}$$
is surjective.
Clearly, the kernel of $\psi$ is isomorphic to 
$\pitop(\Xaffgal-\fgal^{-1}(D))$.

Over a smooth point of $D$ there are $n-1$ points of $X$.
One of them belongs to the ramification divisor of $f$
and has ramification index equal to $2$, whereas the other
points are unramified, cf. Definition \ref{genericdef}.
It follows that the $\Gamma_i$'s map to 
transpositions in $\Sym_n$.
In particular, the subgroup normally
generated by the $\Gamma_i{}^2$'s inside $\pitop(\Aff^2-D)$
lies in $\ker\psi$.

The following lemma extends a result of \cite[Chapter 0.3]{mote}.
\begin{Lemma}
  \label{pioneisoms}
  Let $G_n$ be the subgroup normally generated by
  the $\Gamma_i{}^2$'s in $\pitop(\Aff^2-D)$.
  We view $\Sym_{n-1}$ as a subgroup of $\Sym_n$ fixing
  one of the letters.
  Define $G_{n-1}$ to be the following normal subgroup 
  of $\psi^{-1}(\Sym_{n-1})$ 
  $$
  G_{n-1} \,:=\, \ll G_n,\, \gamma\Gamma_i\gamma^{-1}\,|\,
  \psi(\gamma\Gamma_i\gamma^{-1})\in\Sym_{n-1},\,
  \gamma\in\pitop(\Aff^2-D) \gg .
  $$ 
  Then there exist isomorphisms
  $$\begin{array}{lcl}
    \blowup
    \ker\psi\,/ \,G_n &\iso& \pitop(\Xaffgal)\\
    \psi^{-1}(\Sym_{n-1})\,/ \,G_{n-1} &\iso& \pitop(\Xaff).
  \end{array}$$
\end{Lemma}
\prf
It is shown in \cite[Chapter 0.3]{mote} that the
$\Gamma_i{}^2$ lift to loops around the ramification
divisor $\fgal^{-1}(D)$  and that these loops
generate the kernel of the surjective homomorphism
from $\ker\psi\iso\pitop(\Xaffgal-\fgal^{-1}(D))$ onto
$\pitop(\Xaffgal)$.

The same arguments prove that $G_{n-1}$ is the kernel of 
the surjective homomorphism from 
$\psi^{-1}(\Sym_{n-1})\iso\pitop(\Xaff-f^{-1}(D))$ onto
$\pitop(\Xaff)$.
\qed\medskip

The isomorphism $\ker\psi/ G_n\iso\pitop(\Xaffgal)$ shows
that there is a short exact sequence
\begin{equation}
\label{moteeq}
\begin{array}{ccccccccc}
  1&\to&\pitop(\Xaffgal)&\to&\pitop(\Aff^2-D,\,x_0)/\ll\Gamma_i{}^2\gg
  &\stackrel{\psi}{\to}&\Sym_n&\to&1.
\end{array}
\end{equation}
It is the starting point of the computations 
in \cite[Chapter 0.3]{mote}.
In fact, they prove that this short exact sequence splits.
After that, they use a Reidemeister-Schreier rewriting process
to obtain a presentation of $\pitop(\Xaffgal)$.

\begin{Definition}
  \label{caffinedef}
  We define $\Caffine$
  to be the subgroup normally generated by the following 
  elements inside $\pitop(\Aff^2-D,x_0)$
  $$\begin{array}{ll}
    \left[\gamma\Gamma_i\gamma,\,\Gamma_j{}^{-1}\right]&
    \mbox{if $\psi(\gamma\Gamma_i\gamma^{-1})$ and $\psi(\Gamma_j)$
    are disjoint transpositions}\\
    \langle\gamma\Gamma_i\gamma,\,\Gamma_j{}^{-1}\rangle&
    \mbox{if $\psi(\gamma\Gamma_i\gamma^{-1})$ and $\psi(\Gamma_j)$
    have precisely one letter in common}
  \end{array}$$
  where $\gamma$ runs through $\pitop(\Aff^2-D,x_0)$.
\end{Definition}

\begin{Question}
\label{mainquestion}
  Is the group $\Caffine$ trivial for every good generic projection?
\end{Question}

This is the main problem that we cannot solve.
However, we will see in Section \ref{examplesection} that the answer 
to this question is positive in all known examples up to now.
So the reader should bear in mind that we hope for $\Caffine$
to be trivial.
\medskip

Since $\Caffine$ lies inside $\ker\psi$, the short exact sequence 
(\ref{moteeq}) induces a short exact sequence
\begin{equation}
\label{neweq}
\begin{array}{ccccccccc}
  1&\to&\pitop(\Xaffgal)/\Caffine&\to&\displaystyle
  \frac{\pitop(\Aff^2-D,\,x_0)}{\ll\Gamma_i{}^2,\Caffine\gg}
  &\stackrel{\psi}{\to}&\Sym_n&\to&1.
\end{array}
\end{equation}
We will see in the proof of Theorem \ref{affinetheorem} that this
sequence splits.
Moreover, this short exact sequence is completely determined by
$\pitop(\Xaff)$ and $n$.
However, first we need
\begin{Lemma}
  \label{stillonto}
  The surjective homomorphism of Lemma \ref{pioneisoms} 
  from $\psi^{-1}(\Sym_{n-1})$ onto $\pitop(\Xaff)$ factors
  over the quotient by $\Caffine$.
\end{Lemma}
\prf
Without loss of generality we may assume that $\Sym_{n-1}$ is
the subgroup of $\Sym_n$ fixing the first letter.
We have to show that $\Caffine$ is contained in $G_{n-1}$:
If one of the transpositions $\psi(\gamma\Gamma_i\gamma^{-1})$ or 
$\psi(\Gamma_j)$ lies in $\Sym_{n-1}$ this is true.

Hence we may assume that both transpositions do not lie in $\Sym_{n-1}$.
But this means that they have precisely the first letter in common.
We write $a:=\gamma\Gamma_i\gamma^{-1}$ and $b:=\Gamma_j$.
We have to show that $\langle a,b\rangle$ lies in $G_{n-1}$.
Modulo $G_n$ the elements $a$ and $b$ equal their inverses and so
$\langle a,b\rangle$ is congruent to 
$(aba^{-1})\cdot(b^{-1}a^{-1}b)$.
Both factors lie in $G_{n-1}$ and so the same is true for the product.
\qed\medskip

We now come to our first main result.

\begin{Theorem}
  \label{affinetheorem}
  For a good generic projection $f:X\to\PP^2$ of degree $n$ there exists
  an isomorphism
  $$\begin{array}{ccc}
    \pitop(\Xaffgal)/\Caffine &\iso& \myBigKernel(\pitop(\Xaff),\,n) .
  \end{array}$$
\end{Theorem}
\prf
We consider the short exact sequence (\ref{neweq}).
Let us denote the group in the middle by $A$.
It is generated by the $\Gamma_i$'s, $i=1,...,\deg D$.
We have already seen that the $\Gamma_i$'s and hence
also their conjugates map to  transpositions 
under $\psi$.

First, we split (\ref{neweq}):
For each transposition $\tau$ in $\Sym_n$ we choose a
conjugate $s(\tau)$ of some $\Gamma_i$ in $A$ that maps
to this transposition under $\psi$, i.e.,
$\psi(s(\tau))=\tau$.
Since we quotiented by $\Caffine$ this map $s$, which
is defined on transpositions, extends to a homomorphism
$$\begin{array}{ccccc}
s &:&\Sym_n&\to& A\,\stackrel{{\rm def}}{=}\,
\pitop(\Aff^2-D,\,x_0)/\ll\Gamma_i{}^2,\Caffine\gg
\end{array}$$
that splits $\psi$.

Next, we want to show that $A$ is a quotient of $\EnSym_n(d)$:
We choose for every $\Gamma_i$ an element
$\sigma_i$ in $\Sym_n$ such that 
$\sigma_i\psi(\Gamma_i)\sigma_i{}^{-1}=(1\,2)$.
We define
$$\begin{array}{cccc}
  s_i &:=& s(\sigma_i)\cdot\Gamma_i\cdot s(\sigma_i)^{-1} 
  &i=1,...,\deg D.
\end{array}$$
Clearly, $A$ is generated by the $s_i$'s and $s(\Sym_n)$.
Recall Definition \ref{ensymdef}, where we
defined the group $\EnSym_n(d)$ in terms of 
$d$ generators $s_i$ and the group $\Sym_n$ and 
relations between them.
We define a map
$$\begin{array}{cccl}
  \EnSym_n(d) &\to& A\\
  s_i&\mapsto &s_i&i=1,...,d\\
  \sigma&\mapsto&s(\sigma)&\forall\sigma\in\Sym_n
\end{array}$$
By definition of $\Caffine$, this map extends to a 
surjective homomorphism of groups.

We let $F_d$ with $d=\deg D$ be the free group with basis
$\Gamma_i$, $i=1,...,d$.
The group $\pitop(\Aff^2-D,x_0)$ is a quotient of $F_d$ by
a normal subgroup $R$, cf. Lemma \ref{vankampenlemma}.
We factor this presentation via
$$
\begin{array}{cccccl}
  F_d\ast\Sym_n &\stackrel{\phi}{\to}& \EnSym_n(d) &\to & A\\
  \Gamma_i &\mapsto&s(\sigma_i)^{-1}s_i s(\sigma_i)&\mapsto&\Gamma_i
  &i=1,...,d\\
  \sigma&\mapsto&\sigma&\mapsto&s(\sigma)& \forall \sigma\in\Sym_n
\end{array}$$
By Theorem \ref{rtv}, the group
$\EnSym_n(d)$ is isomorphic to $\myExt(F_{d-1},n)$, where
$F_{d-1}$ is the free group of rank $d-1$.
We will tacitly use this isomorphism. 
In particular, it makes sense to speak of affine subgroups
of $\EnSym_n(d)$, cf. Definition \ref{affinesubdef}.

The group $A$ is a semidirect product of $K:=\pitop(\Xaff)/\Caffine$ 
by $\Sym_n$.
Let $\varphi$ be the action of $\Sym_n$ on the kernel $K$.
For $k\in K$ and $\sigma\in\Sym_n$, we write $[k,\sigma]$ 
for $k\cdot\varphi(\sigma)(k^{-1})$.
Let $G_{n-1}'$ be the image of $G_{n-1}$ in $A$, cf.
Lemma \ref{pioneisoms}.
It is easy to see that $[K,\Sym_{n-1}]$ is equal to the 
intersection $G_{n-1}'\cap K$.
By Lemma \ref{pioneisoms} and Lemma \ref{stillonto}, the
quotient $\psi^{-1}(\Sym_{n-1})/G_{n-1}'$ is isomorphic
to $\pitop(\Xaff)$.
Hence, 
\begin{equation}
\label{symneq}
\begin{array}{ccccc}
  K\,/\,[K,\Sym_{n-1}] &\iso& 
  \psi^{-1}(\Sym_{n-1})\,/\,G_{n-1}'
  &\iso&\pitop(\Xaff),
\end{array}
\end{equation}
where the first isomorphism follows from the second homomorphism
theorem.

Consider the normal subgroup $R':=\phi(R)$ of $\EnSym_n(d)$.
We want to show that $R'$ is an affine subgroup of $\EnSym_n(d)$.
By Lemma \ref{vankampenlemma}, the group $R$ is the normal closure
of elements of a very special type.
It is easy to check that these special elements 
map to elements of the form 
$(1,...,r,1,...,r^{-1},1,...)$ under $\phi$ and the isomorphism
of Theorem \ref{rtv}.
Clearly, such elements generate $R'$.
If we choose another embedding of $\Sym_{n-1}$ into $\Sym_n$,
then formula (\ref{symneq}) still holds true.
This implies that $R'$ is $\Sym_n$-invariant. 
Hence, $R'$ is an affine subgroup.

By Theorem \ref{affinesubthm}, the group $K$ is isomorphic
to $\myBigKernel(G,n)$ for a group $G$.
By Corollary \ref{recover}, this group $G$ is
isomorphic to $K/[ K,\Sym_{n-1}]$.
Applying (\ref{symneq}), we obtain 
$$\begin{array}{ccc}
  \pitop(\Xaffgal)/\Caffine &\iso& \myBigKernel(\pitop(\Xaff),\,n).
\end{array}$$
This was to be shown.
\qed

\begin{Corollary}
  \label{affinecor}
  Let $\overline{\Caffine}$ be the image of $\Caffine$ in
  the Abelianisation of $\pitop(\Xaffgal)$, i.e., in 
  $H_1(\Xaffgal,\ZZ)$.
  Then there exists a non-canonical isomorphism
  $$\begin{array}{ccc}
    H_1(\Xaffgal,\,\ZZ)/\overline{\Caffine} &\iso& H_1(\Xaff,\,\ZZ)^{n-1}.
  \end{array}$$
\end{Corollary}
\prf
This follows from Corollary \ref{bigabelianise} applied to
Theorem \ref{affinetheorem}.
\qed

\section{The fundamental group of the Galois closure}
\label{Mainprojective}

We keep all notations introduced so far.
We assume furthermore that the $\Gamma_i$'s form a good ordered system of 
generators, i.e., that $\delta:=\Gamma_1\cdot...\cdot\Gamma_d$ with $d=\deg D$ 
is homotopic to a loop around infinity in $\ell$.

Before stating our main results, we have to prove a couple of technical
lemmas.

\begin{Lemma}
  \label{deltaseq}
  The element $\delta$ lies in $\ker\psi$ and the centre of
  $\pitop(\Aff^2-D,x_0)$.
  There are two short exact and central sequences
  $$\begin{array}{lcccccccccl}
    \blowup
    &0&\to&\langle\delta\rangle&\to&\pitop(\Xaffgal)&\to&\pitop(\Xgal)&\to&1\\
    &0&\to&\langle\bar{\delta}\rangle&\to&\pitop(\Xaff)&\to&\pitop(X)&\to&1&.
  \end{array}$$
  Here, we denote the image of $\delta$ under
  $\ker\psi\to\pitop(\Xaffgal)$ again by $\delta$
  and its image under $\psi^{-1}(\Sym_{n-1})\to\pitop(\Xaff)$
  by $\bar{\delta}$.
\end{Lemma}
\prf
By a theorem of Oka, $\delta$ is a central element, cf. 
\cite[Corollary 8.4]{fl}.
Of course, $\delta$ remains central in every subquotient of
$\pitop(\Aff^2-D)$.

The subgroup generated by $\delta$ is normal and generates
the kernel of the map from $\pitop(\Aff^2-D)$ onto $\pitop(\PP^2-D)$.
Chasing through the commutative diagram
$$\begin{array}{ccccccccc}
  1&\to&\pitop(\Xaffgal)&\to&\pitop(\Aff^2-D)/\ll\Gamma_i{}^2\gg&
  \stackrel{\psi}{\to}&
  \Sym_n&\to&1\\
  &&\downarrow&&\downarrow&&||\\
  1&\to&\pitop(\Xgal)&\to&\pitop(\PP^2-D)/\ll\Gamma_i{}^2\gg&
  \stackrel{\psi}{\to}&
  \Sym_n&\to&1
\end{array}$$
we obtain the first exact sequence.

For the second statement we exhibit
$\pitop(\Xaff)$ and $\pitop(X)$ as quotients of
$\psi^{-1}(\Sym_{n-1})$ by $G_{n-1}$, cf. 
Lemma \ref{pioneisoms}.
Then we proceed as before.
\qed

\begin{Lemma}
  \label{deltalemma}
  We consider $\myKernel(\pitop(\Xaff),n)$ as a subgroup of
  $\pitop(\Xaff)^n$.
  Under the map (\ref{mybigkerneleq}) the 
  element $\delta$ maps 
  $$\begin{array}{lccccc}
    \blowup&
    \pitop(\Xaffgal)/\Caffine &\iso& \myBigKernel(\pitop(\Xaff),\,n) 
    &\to&\myKernel(\pitop(\Xaff),\,n)\\
    &\delta&\mapsto&\delta'&\mapsto&(\bar{\delta},...,\bar{\delta})
  \end{array}$$ 
  where $\delta$ and $\bar{\delta}$ are as in Lemma \ref{deltaseq}.
  
  In particular, $\bar{\delta}^n$ is an element of the kernel
  $\ker(\pitop(\Xaff)\,\onto\,\pitop(\Xaff)^\ab)$.
\end{Lemma}
\prf
The group $E:=\myBigKernel(\pitop(\Xaff),n)\rtimes\Sym_n$ 
is a quotient of $\pitop(\Aff^2-D)$ by Theorem \ref{affinetheorem}.
As $\delta$ is central in $\pitop(\Aff^2-D)$,
its image $\delta'$ is a central in $E$.
In particular, $\delta'$ is $\Sym_n$-invariant.
The map from $\myBigKernel(\pitop(\Xaff),n)$ onto
$\myKernel(\pitop(\Xaff),n)$ is $\Sym_n$-equivariant.
Thus, $\delta'$ maps to a $\Sym_n$-invariant element, 
i.e., to an element of the diagonal in $\pitop(\Xaff)^n$.

There are two maps from $\pitop(\Xaffgal)$ onto $\pitop(\Xaff)$:
One is the map from Lemma \ref{deltalemma} composed
with the projection onto the first factor from $\pitop(\Xaff)^n$ onto 
$\pitop(\Xaff)$.
We have to show that $\delta$ maps to $\bar{\delta}$ under this map.
The second one is the map induced from Lemma \ref{pioneisoms}.
Under this second map, the element $\delta$ maps to $\bar{\delta}$.
Going through the construction of the first map in the proof of
Theorem \ref{affinetheorem}, we see that it is constructed
via Corollary \ref{recover}, which is compatible with 
Lemma \ref{pioneisoms}.
Hence, the first and the second map coincide.
In particular, $\delta$ maps to $\bar{\delta}$.

The element $(\bar{\delta},...,\bar{\delta})$ 
lies in $\myKernel(\pitop(\Xaff),n)$.
By definition, $\bar{\delta}^n$ lies
in 
$\ker(\pitop(\Xaff)\,\onto\,\pitop(\Xaff)^\ab)$.
\qed

\begin{Lemma}
 \label{complement}
 Let $X$ be a smooth projective surface and $D$ be a smooth and ample
 divisor on it.
 If $\pitop(X)$ is finite then so is $\pitop(X-D)$.
 
 Moreover, if $X$ is simply connected then $\pitop(X-D)$ is a finite
 cyclic group whose order is equal to the divisibility index of
 $\OO_X(D)$ in $\Pic(X)$.
\end{Lemma}
\prf
Let us first assume that $\pitop(X)$ is trivial.
As an ample divisor, $D$ is connected.
Thus, $D$ is irreducible, being connected and smooth.
It follows that every finite cover of $X$ branched along $D$ 
is a cyclic Galois cover.

A cyclic cover branched along $D$ is given by a line
bundle $\cal F$ together with an isomorphism 
${\cal F}^{\otimes c}\iso\OO_X(D)$.
Hence the maximal finite quotient of $\pitop(X-D)$ is
cyclic of order equal to the divisibility index $d$
of $\OO_X(D)$ in $\Pic(X)$.

By Nori's theorem \cite[Corollary 2.5]{nori}, the
group $\pitop(X-D)$ is finitely generated and Abelian.
Hence this group is cyclic of order $d$.
This proves the second assertion.

If $\pitop(X)$ is finite then the universal cover $\tilde{X}$ of
$X$ is again a smooth projective surface.
The inverse image $\tilde{D}$ of $D$ on $\tilde{X}$ is a smooth
and ample divisor.
Hence $\pitop(\tilde{X}-\tilde{D})$ is a finite group by what we
have just proved.
But this latter group is a group of finite index in $\pitop(X-D)$.
Hence $\pitop(X-D)$ is finite, which proves the first assertion.
\qed\medskip

The inclusion map of topological spaces induces
a surjective homomorphism from $\pitop(\Xaffgal)$ onto
$\pitop(\Xgal)$.
Recall that we defined (Definition \ref{caffinedef}) 
a certain subgroup $\Caffine$ of $\pitop(\Xaffgal)$
and determined the structure of $\pitop(\Xaffgal)/\Caffine$.

\begin{Definition}
  We denote the image of $\Caffine$ in $\pitop(\Xgal)$
  by $\Cprojective$.
\end{Definition}

Clearly, if $\Caffine$ is trivial then the same is true for
$\Cprojective$.
As already noted above, $\Caffine$ is trivial in all known examples.

Before giving the structure of $\pitop(\Xgal)/\Cprojective$, we have
to introduce a little bit of notation.
First we define
$$\begin{array}{ccc}
 K &:=& \ker(\,\pitop(\Xaff)\,\onto\,\pitop(X)\,)
  \end{array}$$
which is a cyclic group.
For a natural number $m\geq1$ we define the homomorphism
$$\begin{array}{ccccc}
  \kappa_m &:& K^m & \to & \pitop(\Xaff)^\ab \\
  &&(k_1,...,k_m)&\mapsto& \sum_{i=1}^m \bar{k}_i
\end{array}$$
We are now ready to state our main result on 
$\pitop(\Xgal)/\Cprojective$.

\begin{Theorem}
  \label{projectivemain}
  There exists a cyclic group $Z$ and a central short exact
  sequence
  \begin{equation}
  \label{projeqn1}
  \begin{array}{ccccccccc}
      0 &\to &\displaystyle
      \frac{H_2(\pitop(\Xaff),\,\ZZ)}{Z}
      &\to& \pitop(\Xgal)/\Cprojective
      &\to& G &\to&1
    \end{array}
  \end{equation}
  where $G$ is a group that is itself a central extension
  $$\begin{array}{ccccccccc}
    0&\to&\ker\kappa_{n-1}
    &\to& G &\to&\myKernel(\pitop(X),\,n)
    &\to&1.
  \end{array}$$
  The kernel in (\ref{projeqn1}) lies inside
  the commutator subgroup of $\myKernel(\pitop(X),n)$.
\end{Theorem}

\prf
By Theorem \ref{affinetheorem}, the group $\pitop(\Xaffgal)/\Caffine$
is isomorphic to $\myBigKernel(\pitop(\Xaff),n)$.
By Theorem \ref{mybigkernelthm}, this group is a central extension
$$\begin{array}{cccccccccc}
 0&\to&H_2(\pitop(\Xaff),\,\ZZ) &\to& \pitop(\Xaffgal)/\Caffine &
 \stackrel{\psi}{\to}&
 \myKernel(\pitop(\Xaff),\,n)&\to&1
\end{array}$$
where the kernel lies inside the commutator subgroup.

To obtain $\pitop(\Xgal)/\Cprojective$, we take the quotient
of $\pitop(\Xaffgal)/\Caffine$ by $\langle\delta\rangle$.
We let $Z$ be the intersection of $\langle\delta\rangle$
with $H_2(\pitop(\Xaff),\ZZ)$.
Thus, $\pitop(\Xgal)/\Cprojective$ is a central extension of
$H_2(\pitop(\Xaff),\ZZ)/Z$ by 
$G:=\myKernel(\pitop(\Xaff),n)/\psi(\delta)$.
This gives already the first exact sequence.

By definition, $K$ is the kernel $\ker(\pitop(\Xaff)\to\pitop(X))$.
Using the definition of $\myKernel(-,n)$ it is not hard to see
that the kernel of $\myKernel(\pitop(\Xaff),n)\to\myKernel(\pitop(X),n)$
is equal to $\ker\kappa_n$.
The image of $\delta$ in $\pitop(X)$ is trivial and so
$\psi(\delta)$ lies in $\ker\kappa_n$.
Moreover, by Lemma \ref{deltalemma}, the image $\psi(\delta)$
generates the subgroup of $\ker\kappa_n$ that is equal to
$\Delta(K)$, where $\Delta:K\to K^n$ is the map that sends
$k$ to $(k,...,k)$.
It is not difficult to see that $\ker\kappa_n/\Delta(K)$ is
isomorphic to $\ker\kappa_{n-1}$.

Thus $G$ is an extension of
$\ker\kappa_n/\psi(\delta)\iso\ker\kappa_{n-1}$ by 
$\myKernel(\pitop(X),n)$.
\qed

\begin{Corollary}
  \label{finitecor}
  If $\pitop(X)$ is finite then so is $\pitop(\Xgal)/\Cprojective$.
\end{Corollary}

\prf
If $\pitop(X)$ is finite then so is $\pitop(\Xaff)$
by Lemma \ref{complement}.
Thus, $\pitop(\Xaff)^n$ and $H_2(\pitop(\Xaff),\ZZ)$ are finite.
In particular, $\myBigKernel(\pitop(\Xaff),n)$ is finite by 
Theorem \ref{mybigkernelthm}.

Hence $\pitop(\Xaffgal)/\Caffine$ is finite by Theorem \ref{affinetheorem}.
Since $\pitop(\Xgal)/\Cprojective$ is a quotient of this latter group it is 
also finite.
\qed

\begin{Corollary}
  \label{nilpotentcor}
  If $\pitop(X)$ is nilpotent of class $c$ then $\pitop(\Xgal)/\Cprojective$
  is nilpotent of class at least $c$ and at most $c+2$.
\end{Corollary}

\prf
If $H:=\pitop(\Xgal)$ is nilpotent of class $c$ then the same is true for
$H^n$.
Since $\myKernel(H,n)$ is a subgroup of $H^n$, also $\myKernel(H,n)$ is nilpotent
of class at most $c$.
By Theorem \ref{projectivemain}, we obtain 
$G:=\pitop(\Xgal)/\Cprojective$ as two successive central
extensions of $\myKernel(H,n)$ by Abelian groups.
Hence $K$ is nilpotent of class at most $c+2$.

There is a surjective map from $G$ onto
$\myKernel(H,n)$, which has a surjective map onto $H$.
In particular, $G$ is of class at least $c$.
\qed

\begin{Proposition}
  \label{projectiveabelian}
  Let $\overline{\Cprojective}$ be the image of
  $\Cprojective$ in the Abelianisation of $\pitop(\Xgal)$,
  i.e., in $H_1(\Xgal,\ZZ)$.
  Then there exists a non-canonical isomorphism
  $$\begin{array}{ccc}
    H_1(\Xgal,\,\ZZ)\,/\,\overline{\Cprojective} &\iso&
    H_1(X,\,\ZZ)\,\oplus\,(H_1(\Xaff,\,\ZZ))^{n-2} .
  \end{array}$$
  The group $H_1(\Xaff,\ZZ)$ is an extension of 
  a finite cyclic group (of order dividing $n$) 
  by $H_1(X,\ZZ)$.
\end{Proposition}
\prf
By Lemma \ref{deltaseq}, $\delta$ generates the kernel
$\ker(\pitop(\Xaffgal) \onto \pitop(\Xgal))$.
It follows that the image of $\delta$ in 
$H_1(\Xaffgal)$ generates the kernel of the induced
surjective homomorphism from $H_1(\Xaffgal)$ onto
$H_1(\Xgal)$.
Similarly, the image of $\bar{\delta}$ in $H_1(\Xaff)$
generates the kernel from $H_1(\Xaff)$ onto $H_1(X)$.
By abuse of notation we will write $\delta$ and $\bar{\delta}$
also for their images in the Abelianised groups.
By Lemma \ref{deltalemma}, the order of
$\bar{\delta}$ (as an element of $H_1(\Xaff)$) 
divides $n$.

By Corollary \ref{affinecor} there exists an embedding
$$\begin{array}{ccccc}
  \imath&:&H_1(\Xaffgal)\,/\,\overline{\Caffine} &\to&
  (H_1(\Xaff))^n
\end{array}$$
with image $\myKernel(H_1(\Xaff),n)$.
We denote by $p$ the projection of $H_1(\Xaff)^n$ onto its last 
$n-1$ factors.
Then $p$ induces an isomorphism of $\myKernel(H_1(\Xaff),n)$ with
$H_1(\Xaff)^{n-1}$, cf. Proposition \ref{mykernelAbelian}.
By Lemma \ref{deltalemma}, the element $\delta$ maps to 
$(\bar{\delta},...,\bar{\delta})$ under $\imath$.
Thus, $p\circ\imath(\delta)$ is equal to 
$(\bar{\delta},...,\bar{\delta})$ in $\pitop(\Xaff)^{n-1}$.

To obtain $H_1(\Xgal)/\overline{\Cprojective}$, we have to 
form the quotient of $H_1(\Xaff)^{n-1}$ by 
$\langle\delta\rangle$.
The quotient of $H_1(\Xaff)$ by $\bar{\delta}$ is 
isomorphic to $H_1(X)$.
The image of $\delta$ in $H_1(\Xaff)^{n-1}$
lies on the diagonal $(\bar{\delta},....,\bar{\delta})$.

If $G$ is an Abelian group, $\Delta:G\to G^{n-1}$ the map that
sends $g$ to $(g,...,g)$ and if $N$ is a subgroup of $G$ then
$G^{n-1}/\Delta(N)$ is isomorphic to $G^{n-2}\oplus G/N$.
The proposition then follows if we apply this fact to
$G:=H_1(\Xaff)$ and $N:=\langle\delta\rangle$.\qed

\begin{Corollary}
  \label{betti}
  The rank of
  $H_1(\Xgal,\ZZ)/\overline{\Cprojective}$
  as an Abelian group is equal to $(n-1)$ times 
  the rank of $H_1(X,\ZZ)$.\qed\medskip
\end{Corollary}

\begin{Corollary}
  \label{notorsion}
  Let $X$ be a surface such that $H_1(X,\ZZ)\iso\ZZ^{b_1}$.
  Let $\cal L$ be a sufficiently ample line bundle on
  $X$ and $f:X\to\PP^2$ be a good generic projection with
  respect to $\cal L$.
  We denote by $n$ the self-intersection number 
  and by $d$ the divisibility index of $\cal L$ in $\Pic(X)$.
  Then there are isomorphisms 
  $$\begin{array}{llcl}
    \blowup&
    H_1(\Xaffgal,\,\ZZ)/\overline{\Caffine} &\iso&
    (\Cycl{d})^{n-1}\oplus \ZZ^{b_1(n-1)}\\
    &H_1(\Xgal,\,\ZZ)/\overline{\Cprojective} &\iso& 
    (\Cycl{d})^{n-2}\oplus \ZZ^{b_1(n-1)} .
  \end{array}$$
  These isomorphisms are not canonical.
\end{Corollary}
\prf
Abelianising the short exact sequence of Lemma \ref{deltaseq} 
we obtain 
$$\begin{array}{ccccccccc}
  0&\to&\phi(Z)&\to&\pitop(\Xaff)^\ab&\to&\pitop(X)^\ab&\to&0
\end{array}$$
where $Z:=\langle\bar{\delta}\rangle$ with $\bar{\delta}$
the loop around infinity as in Lemma \ref{deltaseq}
and where $\phi$ denotes the map from $\pitop(\Xaff)$
onto its Abelianisation.
As $\pitop(X)^\ab\iso\ZZ^{b_1}$ is a free Abelian group 
this sequence splits, i.e.,
$\pitop(\Xaff)^\ab\iso\phi(Z)\oplus\ZZ^{b_1}$.

Hence $\phi(Z)$ is a quotient of $\pitop(\Xaff)$.
As in the proof of Lemma \ref{complement} we conclude that
the order of $\phi(Z)$ is equal to the divisibility index
$d$ of $\cal L$.

Using Corollary \ref{affinecor} and 
Proposition \ref{projectiveabelian} finishes the proof.
\qed

\begin{Proposition}
  \label{simplyconnected}
  Let $X$ be a simply connected surface and $\cal L$ a line bundle that
  defines a good generic projection.
  We let $d$ be the divisibility index of $\cal L$ in $\Pic(X)$.
  We denote the self-intersection of $\cal L$ by $n$.
  Then there exist isomorphisms
  $$\begin{array}{llcl}
    \blowup&
    \pitop(\Xaffgal) / \Caffine   &\iso& (\Cycl{d})^{n-1}\\
    &\pitop(\Xgal) / \Cprojective  &\iso& (\Cycl{d})^{n-2}.
  \end{array}$$
  These isomorphisms are not canonical.
\end{Proposition}
\prf
Let $f:X\to\PP^2$ be a good generic projection corresponding to
$\cal L$.
Clearly, the degree $n=\deg f$ is equal to the self-intersection
of $\cal L$.

The inverse image $H:=f^{-1}(\ell)$ of the line at infinity
is a smooth and ample divisor on $X$ by Bertini's theorem.
By definition of $f$, the line bundle $\OO_X(H)$ is 
isomorphic to $\cal L$.
Since $\Xaff=X-H$, Lemma \ref{complement} tells us that
$\pitop(\Xaff)$ is cyclic of order $d$.

Applying Theorem \ref{affinetheorem} and Corollary \ref{bigcyclic},
we obtain
$$\begin{array}{ccccc}
  \pitop(\Xaffgal)/\Caffine &\iso& \myBigKernel(\pitop(\Xaff),\,n)
  &\iso&(\Cycl{d})^{n-1} .
\end{array}$$
In particular, $\pitop(\Xaffgal)/\Caffine$ is Abelian.
So, the same is true for $\pitop(\Xgal)/\Cprojective$.
Thus we can apply Proposition \ref{projectiveabelian} to 
determine $\pitop(\Xgal)/\Cprojective$.
\qed\medskip

If $X$ is not simply connected, we now determine
$H_2(\pitop(\Xaff),\ZZ)$ up to extension.
This will be useful for the examples in the following
section.
\begin{Lemma}
  \label{spectrallemma}
  There exists a short exact sequence of Abelian groups
  \begin{equation}
  \label{h2ext}
  \begin{array}{cccccccc}
     0&\to&\displaystyle\frac{\pitop(X)^\ab\otimes_\ZZ
      K}{\imath(H_3(\pitop(X),\,\ZZ))}
     &\to&H_2(\pitop(\Xaff),\,\ZZ)\\
     &&&\to&\ker ( H_2(\pitop(X),\ZZ)\,\to\, K)&\to&0.
    \end{array}
  \end{equation}
\end{Lemma}
\prf
The kernel $K$ is cyclic.
Hence $H_i(K,\ZZ)=0$ for $i\neq0,1$.
This implies that the spectral sequence
$$
E_{pq}^2\,:=\,H_p(\pitop(X),\,H_q(K,\ZZ))\,\Rightarrow\, 
H_{p+q}(\pitop(\Xaff),\,\ZZ)
$$
has only two rows.
This yields a long exact sequence
$$\begin{array}{l}
  \blowup...\,\to\,H_3(\pitop(X),\ZZ)\,\stackrel{\imath}{\to}\,
   H_1(\pitop(X),K)\,\to\, 
   H_2(\pitop(\Xaff),\ZZ) \,\to\,
   H_2(\pitop(X),\ZZ)\,\to\\
   \,\to\,H_1(K,\ZZ)\,\to\,H_1(\pitop(\Xaff),\ZZ)\,\to\,...
\end{array}$$
Using the universal coefficient formula we obtain isomorphisms 
$$
H_1(\pitop(X),\,K)\,\iso\, H_1(\pitop(X),\,\ZZ)\otimes_\ZZ K \,\iso\,
\pitop(X)^\ab\otimes_\ZZ K\,.
$$
This yields the short exact sequence (\ref{h2ext}).
\qed\medskip

We end this section by a remark on Teicher's conjecture.
In \cite{te}, Teicher observed that $\pitop(\Aff^2-D)$ is virtually 
solvable in some examples and asked whether this might be a general
phenomenon.

\begin{Proposition}
  \label{teicherconjecture}
  If $\pitop(X)$ is not virtually solvable then neither is
  $\pitop(\Aff^2-D)$.
\end{Proposition}

\prf
This follows from the fact that $\pitop(\Aff^2-D)$ has 
$\myBigKernel(\pitop(\Xaff,n))$ as subquotient.
This latter group has $\myKernel(\pitop(X),n)$ as quotient
which has $\pitop(X)$ as a quotient.
Hence, if $\pitop(\Aff^2-D)$ were virtually solvable then
also $\pitop(X)$ would have to be.
\qed\medskip

The fundamental group of a ruled surface over a curve of genus
$\geq2$ contains a free group of rank $2$.
Hence $\pitop(\Aff^2-D)$ for the generic projection from
a ruled surface over a curve of genus $\geq2$ is not virtually
solvable.

\section{Examples}
\label{examplesection}

In this section we discuss all the examples that have been computed so far.
We will see that Question \ref{mainquestion} has a positive answer in
all cases.
We also show how methods from homological algebra simplify and clarify
the computations.

\subsection*{Minimal rational surfaces}

The fundamental groups of Galois closures of generic projections 
from $\PP^2$, $\PP^1\times\PP^1$ and the Hirzebruch surfaces $\FF_e$
have been calculated in a series of papers by Moishezon, Teicher and
Robb.

\begin{Proposition}
  Let ${\cal L}:=\OO_{\PP^2}(k)$ with $k\geq5$ on $X:=\PP^2$.
  Let $f:X\to\PP^2$ be a good generic projection with respect
  to a generic $3$-dimensional subspace of $H^0(X,{\cal L})$
  and $\Xgal$ be the corresponding Galois closure.
  Then there are isomorphisms
  $$\begin{array}{llclc}
    \blowup&
    \pitop(\Xaffgal)/\Caffine &\iso& (\Cycl{k})^{k^2-1}\\
    &\pitop(\Xgal)/\Cprojective &\iso&(\Cycl{k})^{k^2-2}&.
  \end{array}$$
\end{Proposition}
\prf
For $k\geq5$ the line bundle $\cal L$ is sufficiently ample by
Remark \ref{suffample}.
The morphism $f:X\to\PP^2$ has degree equal to the self-intersection
of $\cal L$, which is $k^2$.
The divisibility index of $\cal L$ in $\Pic(X)$ is equal to $k$.
Then we apply Proposition \ref{simplyconnected}.
\qed

\begin{Remark}
  The computations in \cite{mote2} show that the two groups $\Caffine$ 
  and $\Cprojective$
  are trivial already for $k\geq3$, i.e., Question \ref{mainquestion} has a
  positive answer in these cases.
\end{Remark}

\begin{Remark}  
  \label{veronese}
  For $k=2$, the surface $\Xgal$ is the Galois closure of a generic projection
  from the Veronese surface $V_4$ of degree $4$ in $\PP^5$.
  Moishezon and Teicher proved that $\Xgal$ is an Abelian surface, 
  cf. \cite[Proposition 1]{mote2}.
  In particular, $\pitop(\Xgal)$ is isomorphic to $\ZZ^4$,
  whereas we would only predict $(\Cycl{2})^2$ in this case
  (our results are not valid for projections of such small degree).
  However, this example should not be discouraging in view of 
  Question \ref{mainquestion}. 
  The surface $V_4$ has to be excluded from many statements in
  classical algebraic geometry.
  For example, a generic projection from $V_4$ is the only 
  counter-example to Chisini's conjecture \cite{ku}.
\end{Remark}

\begin{Proposition}
  Let ${\cal L}:=\OO_{\PP^1\times\PP^1}(a,b)$ with $a,b\geq5$ on 
  $X:=\PP^1\times\PP^1$.
  Let $f:X\to\PP^2$ be a good generic projection with respect
  to a generic $3$-dimensional subspace of $H^0(X,{\cal L})$
  and $\Xgal$ be the corresponding Galois closure.
  Then there are isomorphisms
  $$\begin{array}{llclc}
    \blowup&
    \pitop(\Xaffgal)/\Caffine &\iso& (\Cycl{\gcd(a,b)})^{2ab-1}\\
    &\pitop(\Xgal)/\Cprojective &\iso& (\Cycl{\gcd(a,b)})^{2ab-2}&.
  \end{array}$$
\end{Proposition}
\prf
For $a,b\geq5$ the line bundle $\cal L$ is sufficiently ample by
Remark \ref{suffample}.
The morphism $f$ has degree equal to the self-intersection
of $\cal L$, which is $2ab$.
The divisibility index of $\cal L$ in $\Pic(X)$ is equal to $\gcd(a,b)$.
Then we apply Proposition \ref{simplyconnected}.
\qed

\begin{Remark}
  The results of \cite{mote} and \cite{mote4} show that 
  $\Caffine$ and $\Cprojective$ are trivial, i.e., Question \ref{mainquestion}
  has a positive answer in these cases.
\end{Remark}

Let us introduce some notation before computing the next example.
We let $X:=\FF_e:=\PP(\OO_{\PP^1}\oplus\OO_{\PP^1}(-e))$ be the 
$e$.th Hirzebruch surface.
We regard $\FF_e$ as a projectivised $\PP^1$-bundle over $\PP^1$. 
We denote by $F$ the class of a fibre and by $H$ the class of the 
tautological bundle inside $\Pic(\FF_e)$.
For $a\gg0$ and $b\gg0$ the bundle 
${\cal L}:=\OO_{\FF_e}(aH+bF)$ will give rise to good generic
projections.
\begin{Proposition}
  Let ${\cal L}:=\OO_{\FF_e}(aH+bF)$ be a sufficiently ample
  line bundle on $X:=\FF_e$.
  Let $f:X\to\PP^2$ be a good generic projection with respect
  to a generic $3$-dimensional subspace of $H^0(X,{\cal L})$
  and $\Xgal$ be the corresponding Galois closure.
  Then there are isomorphisms
  $$\begin{array}{llclc}
    \blowup&
    \pitop(\Xaffgal)/\Caffine &\iso& (\Cycl{\gcd(a,b)})^{2ab+ea^2-1}\\
    &\pitop(\Xgal)/\Cprojective &\iso& (\Cycl{\gcd(a,b)})^{2ab+ea^2-2}&.
  \end{array}$$
\end{Proposition}
\prf
The degree of $f$ is equal to the self-intersection
of $\cal L$, which is $2ab+ea^2$.
The divisibility index of $\cal L$ in $\Pic(X)$ is equal to $\gcd(a,b)$.
Then we apply Proposition \ref{simplyconnected}.
\qed

\begin{Remark}
  The results of \cite{motero} show that $\Caffine$ and $\Cprojective$
  are trivial, i.e., Question \ref{mainquestion} has a positive
  answer in these cases.
\end{Remark}

\subsection*{Product of a curve with $\PP^1$.}

The fundamental groups of the Galois closures of certain generic
projections from $E\times\PP^1$ where $E$ is an elliptic curve,
have been computed by Amram, Goldberg, Teicher and Vishne.
These groups are particularly interesting for us since this is
the first time homological contributions enter the picture.

We let $C$ be a curve of genus $g\geq1$ and 
$X:=C\times\PP^1$.
We let $\cal F$ be a very ample line bundle of degree $k$ on $C$.
The line bundle ${\cal L}:={\cal F}\boxtimes\OO_{\PP^1}(d)$ on $X$
is very ample.
If we choose $k$ and $d$ sufficiently large 
$\cal L$ gives
rise to good generic projections $f:X\to\PP^2$.
The degree $n=\deg f$ is equal to the self-intersection
of $\cal L$, which is equal to $2dk$.

To apply our results, we have to compute $\pitop(\Xaff)$
and $H_2(\pitop(\Xaff),\ZZ)$.

\begin{Lemma}
  The group $\pitop(\Xaff)$ is an extension of
  $\Cycl{d}$ by $\pitop(X)\iso\pitop(C)$.
  If $d$ divides $k$ then this
  extension splits.
  
  Moreover, there is an isomorphism
  $\pitop(\Xaff)^\ab\,\iso\,(\Cycl{\gcd(k,d)})\oplus\ZZ^{2g}$.
\end{Lemma}
\prf
We only do the case $g=1$.
For $g\geq2$ we only have to replace $\CC^2$ by the 
upper half-plane $\HH$.

We let $H$ be a smooth section of $\cal L$.
Then there is a short exact sequence
\begin{equation}
\label{productseq}
\begin{array}{ccccccccc}
0&\to&K&\to&\pitop(\Xaff)&\to&\pitop(X)&\to&1
\end{array}
\end{equation}
where $K$ is a cyclic group in the centre of $\pitop(\Xaff)$,
cf. Lemma \ref{deltaseq}. 

We let $\tilde{X}$ be the universal cover of $X$, i.e.,
$\tilde{X}$ is isomorphic to $\CC^2\times\PP^1$.
The pull-back of $\cal L$ to $\tilde{X}$ is isomorphic
to the line bundle 
$\tilde{\cal{L}}:=\OO_{\PP^1}(d)\boxtimes\OO_{\CC^2}$
as every line bundle on $\CC^2$ is trivial.

The group $C$ corresponds to covers of $\tilde{X}$ ramified
along the inverse image of $H$ on $\tilde{X}$.
Arguing as in the proof of Lemma \ref{complement} we see
that the order of $K$ is equal to the divisibility index
of $\tilde{\cal{L}}$ on $\tilde{X}$, which is equal to $d$.
Hence $K$ is cyclic of order $d$.

If $k$ is divisible by $d$, we can define a 
$d$-fold cover of $X$ by choosing a line bundle 
$\cal F$ such that ${\cal F}^{\otimes d}\iso\cal L$.
Such a cover corresponds to a surjective homomorphism
from $\pitop(X-D)$ onto $K$ compatible with 
(\ref{productseq}).
So in this case the map from $K$ to $\pitop(\Xaff)$ 
is a split injection. 

Abelianising (\ref{productseq}) we obtain $\pitop(\Xaff)^\ab$
as an extension of a finite cyclic group $K'$ by $\ZZ^{2g}$.
This extensions splits since the quotient is free Abelian.
Hence $K'$ is a quotient of $\pitop(\Xaff)^\ab$.
The order of $K'$ is given by the divisibility index of
$\cal L$, which is equal to $\gcd(k,d)$.
\qed

\begin{Lemma}
  There is an isomorphism
  $H_2(\pitop(\Xaff),\ZZ)\,\iso\,(\Cycl{d})^{2g}\oplus\ZZ$.
\end{Lemma}
\prf
The curve $C$ has genus $g\geq1$ and so it is an 
Eilenberg-MacLane space.
Hence the group homology $H_i(\pitop(C),\ZZ)$ is isomorphic 
to the singular homology $H_i(C,\ZZ)$.
In particular, the group $H_3(\pitop(X),\ZZ)$ vanishes and
$H_2(\pitop(X),\ZZ)\iso\ZZ$.
Moreover, $H_1(\pitop(X),\ZZ)\iso\pitop(X)^\ab\iso\ZZ^{2g}$.
By the previous lemma, the kernel $K$ of 
$\pitop(\Xaff)\to\pitop(X)$ is isomorphic to $\Cycl{d}$.

We plug this data into (\ref{h2ext})
and see that $H_2(\pitop(\Xaff),\ZZ)$ is an extension of
$\pitop(X)^\ab\otimes \Cycl{d}$ by 
$L:=\ker(H_2(\pitop(X),\ZZ)\to \Cycl{d})$.
Since $H_2(\pitop(X),\ZZ)\iso \ZZ$, the group $L$
is infinite cyclic, i.e., abstractly isomorphic to $\ZZ$.
So we obtain a short exact sequence of Abelian groups
$$\begin{array}{cccccccccc}
  0&\to&(\Cycl{d})^{2g}&\to&H_2(\pitop(\Xaff),\,\ZZ)&\to&\ZZ&\to&0
\end{array}$$
which splits since the quotient is free Abelian.
\qed\medskip

For a natural number $n\geq1$ we define the homomorphism
$$\begin{array}{ccccc}
  \kappa_n &:& (\Cycl{d})^n &\to& \Cycl{\gcd(d,k)}\\
  &&(x_1,...,x_n) &\mapsto& \sum_{i=1}^n \bar{x}_i 
\end{array}$$
The following statement is now an easy application
of Theorem \ref{projectivemain} and 
Proposition \ref{projectiveabelian}.

\begin{Proposition}
  For a generic projection from $X=C\times\PP^1$ there
  exists a cyclic group $Z$ such that we obtain
  $\pitop(\Xgal)/\Cprojective$ as an extension
  \begin{equation}
   \label{curvetimesp1} 
   \begin{array}{cccccccccc}
    0&\to&\left((\Cycl{d})^{2g}\oplus\ZZ\right)/Z
    &\to&\pitop(\Xgal)/\Cprojective &\to&
    G&\to&1
    \end{array}
   \end{equation}
  where $G$ is a central extension of the form
  $$\begin{array}{cccccccccc}
    0&\to&\ker\kappa_{2dk-1}&\to&G &\to&
    \myKernel(\pitop(C),\,2dk)&\to&1.
  \end{array}$$
  Abelianising we obtain an isomorphism
  $$\begin{array}{ccc}
     H_1(\Xgal,\,\ZZ)/\overline{\Cprojective} &\iso&
     (\Cycl{\gcd(d,k)})^{2dk-2}\,\oplus\,\ZZ^{2g(2dk-1)}.
  \end{array}\qed$$
\end{Proposition}

\begin{Remark}
  The results of \cite{ag} and \cite{agtv} show that $\Caffine$
  and $\Cprojective$ are trivial in the case $g=1$ and $d=1$, i.e.,
  Question \ref{mainquestion} has again a positive answer. 
  Also, in their computations it turns out that $Z$ kills the
  kernel of (\ref{curvetimesp1}), i.e.,
  $\pitop(\Xgal)$ is then isomorphic to
  $\myKernel(\pitop(X),2k)\iso\ZZ^{4k-2}$. 
\end{Remark}

\end{document}